
\def\ph{\varphi}
\def\fx{\Sigma_X}
\def\fy{\Sigma_Y}
\def\Z{{\bf Z}}
\def\P{{\bf P}}
\def\f{\Sigma}
\def\Q{{\bf Q}}
\def\RelInt{\mathop{\rm Rel\,Int}\nolimits}
\def\NE{\mathop{\rm NE}\nolimits}
\def\N{{\cal N}_1}
\def\A{{\cal A}_1}
\def\PC{{\rm PC}}

\hsize=15cm
\vsize=20.5cm
\hoffset=0cm
\parskip 5pt plus 1pt

\catcode`\|=13 
\font\bff=cmbx10  at 14pt
\font\bbk=cmssi10 at 12pt
\font\eightrm=cmr8
\font\eighti=cmmi8
\font\eightsy=cmsy8
\font\eightbf=cmbx8
\font\eighttt=cmtt8
\font\eightit=cmti8
\font\eightsl=cmsl8
\font\sixrm=cmr6
\font\sixi=cmmi6
\font\sixsy=cmsy6
\font\sixbf=cmbx6
\skewchar\eighti='177 \skewchar\sixi='177
\skewchar\eightsy='60 \skewchar\sixsy='60

\def\eightpoint{%
  \textfont0=\eightrm \scriptfont0=\sixrm \scriptscriptfont0=\fiverm
  \def\rm{\fam0\eightrm}%
  \textfont1=\eighti \scriptfont1=\sixi \scriptscriptfont1=\fivei
  \def\oldstyle{\fam1\eighti}%
  \textfont2=\eightsy \scriptfont2=\sixsy \scriptscriptfont2=\fivesy
  \textfont\slfam=\eightit
  \def\sl{\fam\itfam\eightit}%
  \textfont\slfam=\eightsl
  \def\sl{\fam\slfam\eightsl}%
  \textfont\bffam=\eightbf \scriptfont\bffam=\sixbf
  \scriptscriptfont\bffam=\fivebf
  \def\bf{\fam\bffam\eightbf}%
  \textfont\ttfam=\eighttt
  \def\tt{\fam\ttfam\eighttt}%
  \font\petcap=cmcsc8
  \abovedisplayskip=9pt plus 2pt minus 6pt
  \abovedisplayshortskip=0pt plus 2pt
  \belowdisplayskip=9pt plus 2pt minus 6pt
  \belowdisplayshortskip=5pt plus 2pt minus 3pt
  \smallskipamount=2pt plus 1pt minus 1pt
  \medskipamount=4pt plus 2pt minus 1pt
  \bigskipamount=9pt plus 3pt minus 3pt
  \normalbaselineskip=9pt
  \setbox\strutbox=\hbox{\vrule height7pt depth2pt width0pt}%
  \normalbaselines\rm}

\font\petcap=cmcsc10

\newskip\afterskip
\catcode`\@=11
\def\p@int{.\par\vskip\afterskip\penalty100} 

\def\p@intir{\discretionary{.}{}{.\kern.35em---\kern.7em}}
\def\pointir{\afterassignment\pointir@\global\let\next=}
\def\pointir@{\ifx\next\par\p@int\else\p@intir\fi\egroup\next}
\catcode`\@=12

\def|{\relax\ifmmode\vert\else\findef\fi}
\def\findef{\errhelp{Cette barre verticale ne correspond ni a un \vert
mathematique
                        ni a une fin de definition, le contexte doit vous
indiquer ce qui manque.
                        Si vous vouliez inserer un long tiret, le codage
recommande est ---,
                        dans tous les cas, la barre fautive a ete supprimee.}%
                        \errmessage{Une barre verticale a ete trouvee en
mode texte}}

\newif\iffrance
\def\francais{\francetrue 
\frenchspacing
\emergencystretch3em}
\def\anglais{\francefalse}

\def\section#1|{\noindent\bgroup\bf
                                 \par\penalty -500
                                 \vskip 1.5ex plus 1ex minus .2ex
                                 \skip\afterskip=0.5ex plus .2ex
                                  #1\pointir}

\def\resume#1|{\penalty 100{\parindent=0pt\baselineskip=9pt
             {\parshape 2 0pt  12.6cm 1.8cm 10.8cm\bf \skip\afterskip=0pt
             \iffrance R\'esum\'e \else Abstract \fi.\hskip13pt{\ept #1 }\par}}
               \penalty -100}

\def\th#1|{\bgroup \sl \def\findef{\egroup\par}
                        \bgroup\petcap
                         \par\vskip 1ex plus 0.5ex minus .2ex
                         \skip\afterskip=0pt
                           #1\pointir}

\def\rque#1|{\bgroup \sl
                     \par\vskip 1ex plus 0.5ex minus .2ex\skip\afterskip=0pt
                          #1\pointir}

\def\dem{\bgroup \sl
                  \par\vskip 1ex plus 0.5ex minus .2ex\skip\afterskip=0pt
                  \iffrance D\'emonstration\else Proof\fi\pointir}

\def\preuve{\bgroup \sl
                  \par\vskip 1ex plus 0.5ex minus .2ex\skip\afterskip=0pt
                  \iffrance Preuve\else Proof\fi\pointir}

\def\qed{\quad\hbox{\hskip 1pt\vrule width 4pt height 6pt
          depth 1.5pt\hskip 1pt}}
\def\findem{\penalty 500 \hbox{\qed}\par\vskip 3pt}

\def\build#1#2\fin{\mathrel{\mathop{\kern0pt#1}\limits#2}}

\def\ept{\eightpoint}

\def\la{\longrightarrow}

\def\article#1|#2|#3|#4|#5|#6|#7|%
    {\leftskip=7mm\noindent
     \hangindent=2mm\hangafter=1
\llap{[#1]\hskip1.35em}{\bf #2}, {\bf #6}. {#3}, {\sl #4},
     {#5}, #7.\par} 
\def\livre#1|#2|#3|#4|#5|%
    {\leftskip=7mm\noindent
    \hangindent=2mm\hangafter=1
\llap{[#1]\hskip1.35em}{\bf #2}, {\bf #5}. {\sl #3}, {#4}.\par
}
\def\div#1|#2|#3|#4|
{\leftskip=7mm\noindent
\hangindent=2mm\hangafter=1
\llap{[#1]\hskip1.35em}{\bf #2}, {\bf#4}. {#3}.\par
}

\anglais
\newif\ifpagetitre
\newtoks\auteurcourant \auteurcourant={\hfil}
\newtoks\titrecourant \titrecourant={\hfil}

\newcount\pagedebut
\newcount\pagefin
\nopagenumbers
\headline{\ifnum\pageno>\pagedebut
            {\ifodd\pageno \vbox to 0pt{}\eightrm
                    \hfil \the\titrecourant\hfil\folio 
             \else \vbox to 0pt {}\eightrm\folio
                    \hfil\the\auteurcourant\hfil\fi}
          \else\hfil\fi}

\long\def\debutexpose{\advance\pageno by\pagedebut
\advance\pageno by -1
\parindent=0pt
\parindent=3em
\tenpoint
\vskip 2.7cm} 

\pagedebut=1

\footnote{}{$^1$~Partially supported by the European project EAGER.}

\auteurcourant={C.~CASAGRANDE}
\titrecourant={On the birational geometry of toric fano 4-folds}

\vskip1.5cm

{\bff On the birational geometry of toric Fano 4-folds}
\vskip0.3cm

\line{\vtop{\hsize=6.6cm{\bf Cinzia CASAGRANDE~$^1$}
\vskip0.3cm
{\ept\obeylines\parskip0pt
Universit\`a di Roma ``La Sapienza'' 
Dipartimento di Matematica 
piazzale Aldo Moro, 2
00185 Rome (Italie)
E-mail : ccasagra@mat.uniroma1.it\par}}
}
\vskip0.8cm

\line{$\kern1.8cm$\hrulefill$\kern1cm$}

\resume{\eightrm \noindent In this Note, we announce a factorization result
for equivariant birational morphisms between toric 4-folds whose
source is Fano: such a morphism is always a composite of 
blow-ups along smooth invariant centers. 
Moreover, we show with a counterexample that, differently
from the 3-dimensional case, even if both
source and target are Fano, the intermediate varieties can not 
be chosen Fano.}|

\vskip0.5cm
\francais

$\kern30pt${\bbk 
Sur la g{\'e}om{\'e}trie
birationnelle des vari{\'e}t{\'e}s 
toriques de Fano de dimension 4}
\vskip0.3cm

\resume{\eightsl \noindent Dans cette Note, nous annon\c{c}ons un r\'esultat de
factorisation pour les morphismes birationnels \'equivariants entre 
vari{\'e}t{\'e}s toriques de dimension 4, ayant comme source une
vari{\'e}t{\'e} de Fano: un tel morphisme est
toujours donn{\'e} par une suite d'{\'e}clatements le long de
sous-vari\'et\'es lisses invariantes. Nous montrons \`a l'aide d'un
contrexemple que, \`a la diff\'erence du cas de la dimension 3, m{\^e}me si la
source et le but
sont de Fano, les vari{\'e}t{\'e}s interm\'ediaires ne
peuvent pas \^etre choisies de Fano.}| 
 
\line{$\kern1.8cm$\hrulefill$\kern1cm$}
\vskip0.5cm


\parindent=0.4cm

\noindent {\bf Version fran\c{c}aise abr\'eg\'ee}

Dans cette Note, on s'int\'eresse aux morphismes birationnels
\'equivariants entre vari\'et\'es toriques de dimension 4, dans le cas
o\`u la source est une vari\'et\'e de Fano. En dimension 3 
Sato~[9], gr\^ace \`a l'{\'e}tude des
applications birationnelles entre  vari\'et\'es toriques de Fano, 
obtient d'une nouvelle fa\c{c}on la classification
des vari\'et\'es toriques de Fano de
dimension 3 et compl\`ete la classification en dimension 4.
En particulier il obtient le:
\th Th\'eor\`eme 2.1 {\rm (Sato [9])}|
Soient $Y$ et $X$ 
deux vari\'et\'es toriques de Fano de dimension 3 et
 $f\colon Y
\rightarrow X$ un morphisme birationnel \'equivariant. Alors il existe
une suite 
$$
 Y = X_r \build\la^{\varphi_r}\fin X_{r-1}
  \build\la^{\varphi_{r-1}}\fin
\cdots\cdots 
\build\la^{\varphi_2}\fin X_{1}
 \build\la^{\varphi_1}\fin
X_{0}=X  
$$
telle que pour tout $i=1,\ldots,r$ $\varphi_i$ est un {\'e}clatement 
le long d'une sous-vari\'et\'e invariante lisse,  $X_i$ est une
vari\'et\'e torique de Fano et
$f=\ph_1\circ \ldots \circ\ph_r$. 
|
Rappelons que en g\'en\'eral un morphisme birationnel ne se factorise pas 
comme une suite d'{\'e}clatements le long de sous-vari\'et\'es lisses, 
m\^eme dans le cas
torique. La d\'emonstration du th\'eor\`eme 2.1 est obtenue en deux 
\'etapes:\hfill\break
1) en supposant $Y$ de Fano, on montre que $f\colon Y
\rightarrow X$ admet une factorisation en une suite d'{\'e}clatements
le long de sous-vari\'et\'es lisses invariantes; \hfill\break
2) on montre que si de plus $X$ est de Fano, toutes les vari\'et\'es
interm\'ediaires doivent \^etre de Fano.

En dimension 4 nous montrons que seul 1) reste vrai:
\th Th\'eor\`eme 2.2|Soient $Y$ et $X$ 
deux vari\'et\'es toriques lisses, compl\`etes, de dimension 4, et
$f\colon Y
\rightarrow X$ 
un morphisme birationnel \'equivariant. Si $Y$ est de Fano, alors il
existe une suite
$$
 Y = X_r \build\la^{\varphi_r}\fin  X_{r-1}
  \build\la^{\varphi_{r-1}}\fin
\cdots\cdots 
\build\la^{\varphi_2}\fin X_{1}
 \build\la^{\varphi_1}\fin 
X_{0}=X  
\eqno(\heartsuit)$$
o\`u pour tout $i=1,\ldots,r$ $\ph_i$ est un \'eclatement le long 
d'une sous-vari\'et\'e lisse invariante, 
$X_i$ est une vari\'et\'e torique compl\`ete lisse et
$f=\ph_1\circ\ldots\circ\ph_r$.
|
On montre ce r\'esultat  \`a l'aide du proc\'ed\'e utilis\'e
par Sato en
dimension 3: on fixe un c\^one $\sigma$ de dimension 4  dans $\fx$ et
on \'etudie les subdivisions possibles 
 de $\sigma$  dans $\fy$. Le fait que $Y$ soit de Fano
impose beaucoup de restrictions sur la combinatoire de l'\'ev\'entail:
\th Proposition 2.3|
Soient $Y$ et $X$ 
deux vari\'et\'es toriques lisses, compl\`etes, de dimension 4, et
$f\colon Y
\rightarrow X$ 
un morphisme birationnel \'equivariant. 
On suppose $Y$ de Fano et on fixe un c\^one $\sigma\in\fx$ de 
dimension 4. Alors il y a seulement 17 subdivisions possibles de
$\sigma$ dans $\fy$. De ces 17 subdivisions possibles, 4  
apparaissent uniquement si $X\simeq {\P^4}$.
|
Les d\'emonstrations du th\'eor\`eme 2.2 et de la proposition 2.3
appara\^itront ailleurs.

Le point 2) est en d\'efaut d\`es la dimension 4:
\th Proposition 3.1|
Il existe deux vari\'et\'es toriques de Fano de dimension 4 $Y$ et
$X$, et un morphisme birationnel \'equivariant
 $f\colon Y
\rightarrow X$, 
tels que
$f$ ne se factorise pas comme une suite d'{\'e}clatements comme dans 
$(\heartsuit)$, avec $X_i$ de Fano pour tout $i=1,\ldots, r$.
|
Le morphisme $f$ est la composition de deux \'eclatements
$Y\rightarrow W\rightarrow X$, o\`u $X$ est ${\P^4}$ \'eclat\'e le
long d'une droite. 


\vskip 6pt
\centerline{\hbox to 2cm{\hrulefill}}
\vskip 9pt

\anglais

Toric Fano varieties are classified in dimension less or equal than 4: 
${\P^1}$ is the only Fano curve, and there are 5 toric Fano surfaces,
18
toric Fano 3-folds~[1,10] and 124 toric Fano 4-folds~[3,9].
In [9], Sato shows some interesting factorization properties
of equivariant birational maps and morphisms between low-dimensional
toric Fano varieties. With these results, he obtains a new proof
of the classification in dimension 3 and completes the classification
in dimension 4.
 In particular, he shows that every birational
equivariant morphism between toric Fano 3-folds is a composite of
blow-ups with smooth invariant centers between toric Fano 3-folds (theorem 2.1), and
addresses the question whether similar results hold in higher dimension.
We show, with a counterexample,
that the same result is false in dimension 4, and announce a weaker 
factorization result in dimension 4 (theorem 2.2).

In the first section we recall some notions of toric geometry, in
particular the definitions of primitive collection and primitive relation,
and their link with toric Mori theory. The second section concerns 
factorization of birational equivariant morphisms between toric Fano
varieties, and the third one contains our counterexample.

\section 1. Preliminaries in toric geometry|

\noindent For all the standard results in toric geometry, we refer to~[5] and~[7].

Let $X$ be an $n$-dimensional toric variety: $X$ is described by a 
finite fan $\fx$ in the
vector space $N_{\Q}=N\otimes_{\Z}\Q$, where $N$ is a free abelian group of
rank $n$.
We'll always assume $X$ smooth and complete, so the support of
$\fx$ is the whole space $N_{\Q}$ and every cone in $\fx$ is
generated by a part of a basis of $N$.
We remember that for each $r=0,\ldots,n$ there is a bijection between
the cones of dimension $r$ in $\fx$ and the orbits of codimension $r$
in $X$; we'll denote by $V(\sigma)$ the closure of the orbit
corresponding to $\sigma\in\fx$, and $V(x)=V(\langle x\rangle)$ in case
of 1-dimensional cones $\langle x\rangle\in\fx$.
For each 1-dimensional cone $\rho\in\fx$, let $v_{\rho}\in\rho\cap N$
be its primitive generator, and 
$  G(\fx)=\{v_{\rho}\,|\,\rho\in\fx, \dim\rho=1\}  $
the set of all generators in $\fx$.

\th Primitive collections|
{\rm The language of primitive collections and primitive relations was introduced by~Batyrev~([2, 3]); it is particularly convenient to describe
fans of toric Fano varieties.
}|

\th Definition 1.1|
{\rm A subset $P=\{x_1,\ldots,x_h\}\subseteq G(\fx)$ is a {\it primitive collection}
for $\fx$ if $\langle x_1,\ldots,x_h\rangle
\notin\fx$, but $\langle x_1,\ldots,\check{x}_i,\ldots,x_h\rangle\in\fx$
for each $i=1,\ldots,h$.}|

\noindent We denote by $\PC(\fx)$ the set of all primitive collections 
for $\fx$.  

\th Definition 1.2|{\rm Let $P=\{x_1,\ldots,x_h\}\subseteq G(\fx)$ be a
primitive collection. Since $X$ is complete, the point
$x_1+\cdots+x_h$ is contained in some cone of $\fx$; let
$\sigma_P=\langle y_1,\ldots,y_k \rangle$ be the unique cone in $\fx$
such that
$ x_1+\cdots +x_h\in \RelInt \sigma_P$, the relative interior of $\sigma_P$. 
Then we get a linear relation
$$
x_1+\cdots +x_h-(a_1y_1+\cdots+ a_ky_k)=0
$$
with $a_i$ a positive integer for each $i=1,\ldots,k$. We call this
relation the {\it primitive relation} associated to $P$.
The {\it degree} of $P$ is the integer $\deg P = h-a_1-\cdots-a_k$.}
|

Let $\A(X)$ be the group of algebraic 1-cycles on $X$ modulo numerical 
equivalence, $\N(X)=\A(X)\otimes_{\Z}\Q$, and $\NE(X)\subset \N(X)$ the
cone of Mori, generated by classes of effective curves. The
$\Q\,$-vector space $\N(X)$ has dimension $\rho_X$, the Picard number of $X$. 
We recall the well-known result:

\th Proposition 1.3|
The group $\A(X)$ is canonically isomorphic to the lattice of
integral relations among the elements of $G(\fx)$.
|

\noindent A relation $\sum_{x\in G(\fx)} a_xx=0$ corresponds to a class
in $\A(X)$ whose intersection with $V(x)$ is $a_x$, for all $x\in
G(\fx)$.
By proposition 1.3, for
every primitive collection $P\in\PC(\fx)$, the associated primitive
relation defines a class $r(P)\in\A(X)$. 
Since the canonical class on $X$ is given by $K_X=-\sum_{x\in G(\fx)}
 V(x)$, for every primitive collection $P$ we have 
$ -K_X\cdot r(P)=\deg P. $

\th  Some links with toric Mori theory|
{\rm The cone of effective curves in a complete toric variety has been
studied by Reid
in~[8]: it is a closed, polyhedral cone, generated by classes
of invariant curves.
Moreover, it follows from Kleiman's criterion of ampleness~[6]
that $X$ is projective if and only if the cone $\NE(X)$ is strictly
convex; in this case, its extremal rays are spanned by 
invariant curves.}
|

\th Theorem 1.4 {\rm (Reid [8], Batyrev [2])}|
The cone of effective curves $\NE(X)$ is generated by primitive
relations:
$$ \NE(X)=\sum_{P\in\PC(\fx)}\Q_{\, \geq 0}\, r(P). $$
|

\noindent This gives an important characterization of
toric Fano varieties:

\th Proposition 1.5 {\rm (Batyrev [3])}|
The variety $X$ is Fano if and only if all primitive
collections in $\fx$ have strictly positive degree.
|

\th Blow-ups|{\rm 
By a smooth equivariant blow-up we mean the blow-up of a smooth toric
variety along an invariant, smooth subvariety. The resulting variety
is clearly a smooth toric variety.}|

Let $f\colon Y\rightarrow X$ be a smooth equivariant blow-up, along a 
subvariety $V(\tau)\subset X$, $\tau=\langle x_1,\ldots,x_h\rangle$.
We recall that $\fy$ is a subdivision of $\fx$, therefore
$G(\fx)\subseteq G(\fy)$.
The set $P=\{x_1,\ldots,x_h\}$ is a primitive collection for $\fy$, with relation
$$
r(P)\colon\quad x_1+\cdots + x_h =x.       \eqno(\star)$$
The divisor $V(x)$ is the exceptional divisor in $Y$ and $r(P)$ corresponds to 
the numerical class of a $\P^1$ contained in a fiber of $f$.
In general, $r(P)$ needs not to be extremal in $\NE(Y)$: 

\th  Proposition 1.6 {\rm (Bonavero [4])}|
Suppose $Y$ that is projective: then $r(P)$ generates an extremal ray 
in $\NE(Y)$ if and only if $X$ is projective.|

\section 2. Results on birational morphism between toric Fano varieties|

\noindent In [9], Sato uses primitive collections and primitive
relations to study
equivariant birational maps between toric Fano varieties. In this way, he
shows a strong property of birational morphisms
between toric Fano 3-folds, without using the classification:

\th Theorem 2.1 {\rm (Sato [9])}|
Let $Y$ and $X$ be two toric Fano 3-folds, and $f\colon Y
\rightarrow X$ a birational equivariant morphism. 
Then there exists a sequence of smooth equivariant blow-ups
$$
 Y = X_r \build\la^{\varphi_r}\fin X_{r-1}
  \build\la^{\varphi_{r-1}}\fin
\cdots\cdots 
\build\la^{\varphi_2}\fin X_{1}
 \build\la^{\varphi_1}\fin
X_{0}=X  
$$
such that for each $i=1,\ldots,r$ $X_i$ is a toric Fano 3-fold and
$f=\ph_1\circ \ldots \circ\ph_r$. 
|

We remark that in general it is not possible to factorize a birational
morphism in a sequence of smooth blow-ups; here not only it is
possible, but also all the intermediate varieties are Fano.
In fact, theorem~2.1 is obtained in two steps:\hfill\break
1)~~given a birational equivariant morphism $f\colon Y\rightarrow X$,
if $Y$ is Fano, it is always possible to find a factorization of $f$
in a sequence of smooth, equivariant blow-ups;\hfill\break
2)~~if moreover $X$ is Fano, all the intermediate varieties must be Fano.

In the 4-dimensional case, only the first part of this result stands true:

\th Theorem 2.2|Let $Y$ and $X$ be nonsingular, complete toric varieties of dimension 4,
and $f\colon Y
\rightarrow X$ a birational equivariant morphism. 
Suppose $Y$ is Fano:
then there exists a sequence of smooth equivariant blow-ups
$$
 Y = X_r \build\la^{\varphi_r}\fin  X_{r-1}
  \build\la^{\varphi_{r-1}}\fin
\cdots\cdots 
\build\la^{\varphi_2}\fin X_{1}
 \build\la^{\varphi_1}\fin 
X_{0}=X  
\eqno(\heartsuit)$$
 such that for each $i=1,\ldots,r$ $\ph_i$ is a smooth 
equivariant blow-up, $X_i$ is a nonsingular, complete toric variety
and $f=\ph_1\circ\ldots\circ\ph_r$.
|
The proof of this result uses the same idea of Sato's proof in
the 3-dimensional case: for a fixed cone $\sigma\in\fx$, we study the
possible subdivisions of $\sigma$ in $\fy$. The hypothesis that $Y$ is
Fano gives many combinatorial restrictions on the possible
subdivisions:

\th Proposition 2.3|
Let $Y$ be a toric Fano 4-fold,
$X$ a nonsingular, complete toric 4-fold
and $f\colon Y\rightarrow X$ a birational equivariant morphism. 
We fix a 4-dimensional cone $\sigma\in\fx$. Then there are only 17
possible subdivisions of $\sigma$ in $\fy$, 4 of which can appear only
if $X\simeq {\P^4}$.
|

The proofs of theorem 2.2 and proposition 2.3 will appear elsewhere.

\section  3. The example|

In this section we are going to show the
\th Proposition 3.1|
There exist two nonsingular toric Fano 4-folds
$Y$ and $X$ and a birational equivariant morphism $f\colon Y
\rightarrow X$, such that $f$ doesn't admit a decomposition in smooth
equivariant blow-up as in
$(\heartsuit)$ with $X_i$ a toric Fano 4-fold  for all 
$i=1,\ldots, r$.
|
\dem
Fix two lines $L_1$, $L_2$ in  ${\P^4}$ which intersect in a
point.
If $\{e_1, e_2, e_3, e_4\}$ is a basis of $N$, then we have
$G(\Sigma_{\P^4})=\{e_0,e_1,e_2,e_3,e_4\}$
with the  only  primitive relation
$ e_0+e_1+e_2+e_3+e_4=0.$
We can suppose that $L_1=V(\langle e_1,e_2,e_3\rangle)$, 
$L_2=V(\langle e_2,e_3,e_4\rangle)$;
 let $X$ be the blow-up of  $\P^4$ along $L_1$.

In [9], Sato studies how primitive collections vary under 
smooth equivariant blow-ups and blow-downs. 
Let $f\colon Z\rightarrow W$ be a blow-up: Sato
 gives an algorythm to compute 
explicitly the primitive collections of $\f_Z$ from the ones of 
$\f_W$, and viceversa. We have used  this algorythm  to
obtain the primitive collections of $X$ and of the
varieties we'll introduce in the sequel.

The fan $\Sigma_X$ of $X$ in $N_{\Q}$ has vertices
G($\Sigma_X$)=$\{e_0,\ldots,e_5\}$ and its primitive relations
are $e_1+e_2+e_3=e_5$ and $e_0+e_4+e_5=0$.
These relations are both extremal; 
$\N(X)$ is a vector space of dimension $\rho_X=2$ and
$\NE(X)\subset\N(X)$
is a 2-dimensional cone.
Furthermore, the primitive relations have respectively degree 2 and 3,
so $X$ is a nonsingular toric Fano 4-fold, by proposition~1.5.
In $X$ we have the invariant subvarieties:

\noindent $\widetilde{L_2}=V(\langle e_2,e_3,e_4\rangle)$, the strict transform
  of $L_2$;\hfill\break
$D=V(e_5)\simeq{\P^2}\times{\P^1}$, the exceptional
  divisor;\hfill\break
$S=V(\langle e_4,e_5\rangle)$, the ${\P^2}$ contained in $D$
  passing through the point $p=\widetilde{L_2}\cap D$;\hfill\break
$C_1=V(\langle e_2,e_4,e_5\rangle)$ and $C_2=V(\langle
  e_3,e_4,e_5\rangle)$, the two invariant lines in $S$ passing through $p$.

Now we blow up $X$ along the curve $\widetilde{L_2}$.
We get a toric 4-fold $W$ whose fan
$\Sigma_W$ has vertices G($\Sigma_W$)=$\{e_0,\ldots,e_6\}$
and primitive relations:
$$\matrix{
(E)~~e_2+e_3+e_4  \hfill&= &  e_6  \qquad \hfill  &\qquad e_1+e_2+e_3  \hfill&= &  e_5 \hfill&\cr
(E)~~e_1+e_6      \hfill&= &  e_4+e_5\qquad  \hfill &\qquad e_0+e_4+e_5  \hfill&= &  0\hfill&\cr
(E)~~ e_0+e_5+e_6  \hfill&= &  e_2+e_3\hfill\cr
}$$
Here the letter $E$ stands for extremal; 
 the decomposition of the last two relations are:
$$
r(\{e_1,e_2,e_3\})= r(\{e_1,e_6\})+r(\{e_2,e_3,e_4\}),\quad 
r(\{e_0,e_4,e_5\})= r(\{e_0,e_5,e_6\})+r(\{e_2,e_3,e_4\}).
$$
Therefore $\NE(W)$ is a simplicial cone in
  $\N(W)$, which has dimension $\rho_W=3$. 
The variety $W$ is not Fano: indeed, there is
a primitive collection $\{e_1, e_6\}$ of degree zero.
The invariant curves whose numerical class in $\NE(W)$
corresponds to the relation $e_1+e_6= e_4+e_5$ are exactly 
$\widetilde{C_1}=V(\langle e_2,e_4,e_5\rangle)$ and 
$\widetilde{C_2}=V(\langle
  e_3,e_4,e_5\rangle)$, the strict transforms of $C_1$ and $C_2$, both
  contained in the surface  $\widetilde{S}=V(\langle e_4,e_5\rangle)$, which
  is now a $\P^2$ blown-up in one point, {\it i.e.\/}  the Hirzebruch
  surface ${\bf F}_1$. The curves $\widetilde{C_1}$ and
  $\widetilde{C_2}$ have anticanonical degree zero.

Finally, let $Y$ be the blow-up of $W$ along the surface $\widetilde{S}$. 
The fan $\Sigma_Y$ of $Y$ in $N$
has vertices  G($\Sigma_Y$)=$\{e_0,\ldots,e_7\}$
and primitive relations:
$$\matrix{
(E)~~ e_4+e_5 \hfill&= &e_7 \hfill &\qquad (E)~~e_2+e_3+e_7 \hfill&= &e_5+e_6\hfill 
&\qquad e_1+e_2+e_3 \hfill&=& e_5 \hfill& \cr
(E)~~ e_1+e_6 \hfill&= &e_7 \hfill &\qquad(E)~~ e_0+e_7 \hfill&= &0\hfill
&\qquad e_2+e_3+e_4 \hfill&=& e_6.\hfill& \cr
}$$
There are 4 extremal classes, and  $\N(Y)$ has dimension $\rho_Y=4$, 
hence the cone $\NE(Y)\subset\N(Y)$ is again simplicial.
The decompositions of the last two relations are:
$$ r(\{e_1,e_2,e_3\})= r(\{e_1,e_6\})+r(\{e_2,e_3,e_7\}),\quad
r(\{e_2,e_3,e_4\})=r(\{e_4,e_5\})+r(\{e_2,e_3,e_7\}).
$$
Since all primitive relations have strictly positive degree,
proposition~1.5 implies that $Y$ is a Fano variety. 
{\it We claim that $Y$ doesn't admit any equivariant blow-down
to a nonsingular toric Fano 4-fold.} 

Indeed, any equivariant blow-down of $Y$ to a smooth
toric 4-fold gives a primitive relation as in $(\star)$; in $\fy$
there are 4 primitive relations of this type.

The two relations
$e_4+e_5=e_7$ and $e_1+e_6=e_7$ give blow-downs respectively
to $W$ and $\overline{W}$, where  $\overline{W}$ 
is obtained from ${\P^4}$ blowing-up 
first $L_2$ and then the strict transform $\widetilde{L_1}$;  
$\overline{W}$ is clearly isomorphic to $W$, therefore it is not Fano.

The other two relations to consider are
$e_2+e_3+e_4=e_6$ and 
$e_1+e_2+e_3=e_5$. 
Suppose, for instance, that $e_1+e_2+e_3=e_5$ comes 
from an equivariant blow-up, as in $(\star)$. 
Then every cone of $\fy$ containing $e_5$ must  
contain also two generators among 
$\{e_1,e_2,e_3\}$.
On the other hand, by Reid's study of
 the geometry of a fan around an extremal wall in~[8], 
the fact that $e_2+e_3+e_7=e_5+e_6$ is extremal
 implies $\langle e_3,e_5,e_6,e_7\rangle\in\fy$;
this gives a contradiction.
The same for  $e_2+e_3+e_4=e_6$.
  
Therefore, the birational morphism $Y\rightarrow X$ given by the
composition of the two blow-ups $Y\rightarrow
W\rightarrow X$ is the example we were looking for.
\findem
\medskip

{\it Aknowledgements.} 
This note was written while I was visiting the Institut Fourier in
Grenoble: I am grateful to this institution for its kind hospitality. 
I would like to thank Lucia Caporaso and
Laurent Bonavero for many valuable conversations and advices, 
 and also for their encouragement.

\bigskip

\centerline{\bf  References}
\medskip

{\baselineskip=2pt\eightrm
\item{[1]}Batyrev, V.V.: Toroidal Fano 3-folds, Math. USSR-Izv. 19,
13--25 (1982).

\item{[2]}Batyrev, V.V.: On the classification of smooth projective
toric varieties, Tohoku Math. J.  43, 569--585 (1991).

\item{[3]}Batyrev, V.V.: On the classification of toric Fano
4-folds, J. Math. Sc. (New York) 94, 1021--1050 (1999). 

\item{[4]}Bonavero L.: Sur des vari{\'e}t{\'e}s toriques non
projectives, Bull. Soc. Math. France 128, 407--431 (2000). 

\item{[5]}Fulton W.: Introduction to Toric Varieties, Ann. of
Math. Studies 131, Princeton Univ. Press (1993).

\item{[6]}Kleiman S.L.:  Towards a numerical theory of ampleness,
Ann. of Math. 84, 293--344 (1966).

\item{[7]}Oda T.: Convex Bodies and Algebraic Geometry - An
Introduction to the Theory of Toric Varieties, Springer-Verlag (1988).

\item{[8]}Reid M.: Decomposition of Toric Morphisms, in Arithmetic
and Geometry, vol. II: Geometry, Prog. Math. 36, 395--418 (1983).

\item{[9]}Sato H.: Toward the classification of higher-dimensional
toric Fano varieties, Tohoju Math. J. 52, 383--413 (2000).

\item{[10]}Watanabe K., Watanabe M.: The classification of Fano
3-folds with torus embeddings, Tokyo J. Math. 5, 37--48 (1982).

}

\end